\documentclass[11pt]{article}
\usepackage{amssymb}
\usepackage{amsmath}
\usepackage{amsthm}
\usepackage[all]{xy}
\usepackage{david}

\newcommand{\Sup}{\mathcal{S}up}
\theoremstyle{definition}
\newtheorem{oppr}[defin]{Open problems}

\begin{document}
\title{Girard couples of quantales}
\author{J.~M.~Egger \and
David Kruml\thanks{Supported by the Ministry of Education
of the Czech Republic under the project MSM143100009.}}
\maketitle

{\bf Keywords:} Girard quantale, quantale, spectrum of
operator algebra, $*$-autonomous functor category.

\begin{abstract}
We introduce the concept of a \emph{Girard couple}, which consists of
two (not necessarily unital) quantales linked by a strong form of
duality.   
The two basic examples of Girard couples arise in the study of
endomorphism quantales and of the spectra of operator algebras.   
We construct, for an arbitrary sup-lattice $S$, a Girard quantale
whose right-sided part is isomorphic to $S$. 
\end{abstract}

\section{Introduction}

Girard quantales were introduced by Yetter to provide semantics for a
certain fragment of {non-commutative linear logic} known as 
\emph{cyclic linear logic} \cite{Gir,Yet}.
They are, essentially, quantales with a {well-behaved} negation
operation, and therefore play a role among all quantales analogous to
that played by complete boolean algebras among frames.  
They are also related to the much older notion of \emph{MV-algebra}
\cite{Cha,Pao}.  

Endomorphisms quantales $\Q(S)$ have been studied by C.~J.~Mulvey and
J.~Wick Pelletier as quantales of \emph{linear relations}
\cite{MuPe92}. 
These quantales admit a \scary{von Neumann duality} between their
right- and left-sided elements; but this can not, in general, be
extended to arbitrary elements of the quantale. 

In this paper, we construct a \emph{predual quantale} $\C(S)$ for the
endomorphism quantale $\Q(S)$ and show that the pairs $(\C(S),\Q(S))$
enjoy properties among all suitable pairs of quantales which are
analogous to those of a Girard quantale. 
In particular, we define a negation operation which extends the
existing von Neumann duality. 

Our construction is reminiscent of one arising in functional analysis,
where the ideal of trace-class operators on a Hilbert space is the
(Banach space) predual of the algebra of all bounded operators on
that Hilbert space. 
By considering appropriate topologies on the preceding algebras, we
can construct further fundamental examples of what we shall call 
\emph{Girard couples}.

We also note that there is a characterisation of Girard couples in
terms of monoidal functors which, in turn, suggests further
generalisations: by considering more complex \scary{gradings} of
quantale structures, or by replacing sup-lattices by objects of an  
other \staut\ category.

\section{Preliminaries}
We review some of the basic definitions and results of quantale theory
which will be extensively used in the sequel.  
Details may be found in \cite{Bar,JoTi84,KaRo97,KrPa07,MuPe01,Rose90}.

The category of complete lattices and supremum-preserving maps will be
denoted $\Sup$; we shall follow the convention of referring to objects
and arrows of $\Sup$ as \emph{sup-lattices} and 
\emph{sup-homomorphisms}, respectively.
The top and bottom elements of a sup-lattice will be denoted $1,0$, 
respectively.
We say that a sup-homomorphism is \emph{strong} if it preserves the
top element. 

The category $\Sup$ has a {\staut\ structure}:  
the \emph{tensor product} of sup-lattices $S$ and $T$, denoted 
$S \otimes T$, is the free sup-lattice with generators  
$\{s \otimes t\mid s\in S, t\in T\}$ satisfying the relations  
\begin{align*}
\bigvee (s_i \otimes t) &= \left(\bigvee s_i\right) \otimes t&
\bigvee (s \otimes t_j) &= s \otimes \left(\bigvee t_j\right)
\end{align*}
for all $s,s_i\in S,t,t_j\in T$;
the \emph{tensor unit} is the two-element chain $\Two=\{0,1\}$;
the \emph{dual} of a sup-lattice $S$ is simply its opposite, 
denoted $S\op$.

We mark elements of $S\op$ with $'$ whenever the distinction from
elements of $S$ is desirable. 
Every sup-homomorphism $f:S\to T$ has a right adjoint $f\adj:T\to S$
which preserves arbitrary infima, and so may be regarded as a
sup-homomorphism $f^*:T\op\to S\op,f^*(x')=f\adj(x)'$; 
this is the \emph{dual} of $f$.

A \emph{quantale} is a sup-lattice $Q$ equipped with an associative
multiplication that distributes joins. 
The right adjoints of $a\cdot\blank$ and $\blank\cdot a$ are denoted  
$\blank\la a,a\ra\blank$, respectively; they can be computed as below.
\begin{align*}
b\la a &=\bigvee\{ c\mid ac\leq b\} &
a\ra b &=\bigvee\{ c\mid ca\leq b\}
\end{align*}
An element $r\in Q$ is said to be \emph{right-sided} if $r1\leq r$.
Similarly, $l\in Q$ is \emph{left-sided} if $1l\leq l$.
A \emph{two-sided} element is both right- and left-sided.
The sets of right-, left-, and two-sided elements are denoted
$\R(Q)$, $\LL(Q)$  and $\T(Q)$, respectively.
Note that the \emph{left annulator} $a\ra 0$ of any $a \in Q$ is
left-sided and that the \emph{right annulator} $0\la a$ is
right-sided. 
Thus the mappings $\blank\ra 0,0\la\blank$ establish a pseudoduality  
between $\R(Q)$ and $\LL(Q)$. 
We write $l\perp r\LRa lr=0$, $r\comp=r\ra 0$, $l\comp=0\la l$. 

We say that $Q$ is: 
\emph{unital} if it has a neutral element, 
\ie, an $e\in Q$, such that $ae=ea=a$ for all $a\in Q$; 
\emph{semiunital} if $r1=r,1l=l$ for all $r\in\R(Q),l\in\LL(Q)$;
\emph{von Neumann} if $\comp$ is a duality between $\R(Q)$ and
$\LL(Q)$. 
In a semiunital quantale $Q$, $1a\geq a$ and $a1\geq a$ hold for every
$a\in Q$. 

A \emph{homomorphism of quantales} is a sup-homomorphism preserving 
the multiplication. 
It is \emph{unital} if it also preserves the neutral element.

A \emph{left $Q$-module} is a sup-lattice $M$ together with an action
of $Q$ on $M$ which respects joins in both variables and satisfies
$(ab)m=a(bm)$ for all $a,b\in Q,m\in M$. 
\emph{Right $Q$-modules} are defined similarly and a
\emph{$Q$-bimodule} is required to also satisfy $(am)b=a(mb)$ for all
$a,b\in Q,m\in M$. 

We say that a left $Q$-module $M$ is:
\emph{unital} if $Q$ is unital and $em=m$ for every $m\in M$; 
\emph{strong} if $1_Qm=1_M$ for every $m\in M,m\neq 0$.
The right adjoints of $a\cdot\blank$ and $\blank\cdot a$ on a
$Q$-bimodule $M$ are also denoted $\blank\la a,a\ra\blank$.
We also write $m\ra n=\bigvee\{ a\in Q\mid am\leq n\}$,
$m\la n=\bigvee\{ a\in Q\mid na\leq m\}$ for $m,n\in M$.

A \emph{homomorphism of left $Q$-modules} $f:M\to N$ is a
sup-homomorphism which satisfies $f(am)=af(m)$ for every 
$a\in Q,m\in M$. 
Homomorphisms of right $Q$-modules and $Q$-bimodules are defined
in a similar way.

An important example is $\Q(S)$, the quantale of endomorphisms of a
fixed sup-lattice $S$ with composition as multiplication and suprema   
calculated pointwise. 
Its right- and left-sided elements are those of the form
\begin{align*}
\rho_x(y)=\begin{cases}
x, & y\neq 0,\\
0, & y=0,
\end{cases}& &
\lambda_x(y)=\begin{cases}
1, & y\nleq x,\\
0, & y\leq x,
\end{cases}
\end{align*}
hence $\T\Q(S)=\Two$; $\Q(S)$ is von Neumann because 
$\rho_x\comp=\lambda_x$, $\lambda_x\comp=\rho_x$.
Moreover $\Q(S)$ is \emph{simple} \cite{PeRo97,Pa97}
and every element $\alpha\in\Q(S)$ can be expressed as below.
$$\alpha=\bigwedge_{x\in S}(\rho_{\alpha(x)}\vee\lambda_x)
=\bigwedge_{x\in S}(\rho_x\vee\lambda_{\alpha\adj(x)})$$

An element $d$ of a quantale $Q$ is called:
\emph{cyclic} if $ab\leq d\LRa ba\leq d$ for all $a,b\in Q$;
\emph{dualizing} if $d\la (a\ra d)=(d\la a)\ra d=a$ for every $a\in Q$. 
$Q$ is called \emph{Girard} if it has a cyclic dualizing element.
In that case we write $a\comp=a\ra d=d\la a$.

The \emph{spectrum} of a C*-algebra $A$, denoted $\Max A$, is the
sup-lattice of all linear subspaces of $A$ which are closed \wrt\ the
norm topology.  
It is a quantale \wrt\ the multiplication 
$ab=\cl\{ AB\mid A\in a,B\in b\}$.
The right- and left-sided elements of $\Max A$ are, respectively, the
closed right and left ideals of $A$. 

Given a von Neumann algebra $M \subseteq \B(H)$, one may consider
either the \emph{weak spectrum} $\Max_wM$ 
or the \emph{ultraweak spectrum} $\Max_{\sigma w}M$; 
the former consists of all linear subspaces of $M$ which are closed
\wrt\ the {weak (operator) topology}, the latter of those which are
closed \wrt\ the somewhat finer {ultraweak topology}.
The weak spectrum of a von Neumann algebra is a von Neumann quantale 
\cite{Pe97}; it follows that the same is true for ultraweak spectra,
which are better suited to our purposes.
[It is well-known that an ideal is ultraweakly closed if and only if
it is weakly closed.]

We recall that a functional $\B(H) \to \CC$ is ultraweakly continuous
if and only if it has the form 
$\sum_{i=1}^\infty x_i(\phi_i,\blank\cdot\psi_i)$
for some orthonormal families $\phi_i,\psi_i\in H$ and coefficients
$x_i\in\CC$ such that $\sum |x_i|<\infty$.
Moreover, a subspace of $\B(H)$ is ultraweakly closed if and only if
it is the intersection of the kernels of some set of ultraweakly 
continuous functionals.

An element $C\in\B(H)$ is said to be \emph{trace-class} if
$$\| C\|_1=\sum_{i=1}^\infty (\phi_i,\sqrt{C^*C}\phi_i)<\infty$$
for some orthonormal basis $\phi_i$ of $H$.
The set of all trace-class elements is an ideal in $\B(H)$ and is
denoted $\C_1(H)$. 
The number $\| C\|_1$ does not depend on the chosen basis and defines
a norm on $\C_1(H)$. 
The $\|\ \|_1$-closed subspaces form a spectrum $\Max_1\C_1(H)$.

The more general \emph{Schatten class} $\C_p(H)$ for $p\geq 1$ is
defined as a set of all elements $C\in\B(H)$ such that 
$$\| C\|_p=\left(\sum_{i=1}^\infty
(\phi_i,\sqrt{C^*C}\phi_i)^p\right)^{1/p}<\infty;$$
subspaces closed in the $\|\ \|_p$-norm form a spectrum,
$\Max_p\C_p(H)$. 

Given a family of Hilbert spaces $H_i$, we can think of the algebra
$\prod\B(H_i)$ as that subalgebra of $\B(\bigoplus H_i)$ consisting of
those operators for which $H_i$ are invariant subspaces;  
similarly, $\bigoplus\C_1(H_i)=\C_1(\bigoplus H_i)\cap\prod\B(H_i)$.
The induced topologies on
$\prod\B(H_i)\subseteq\B(\bigoplus H_i)$,
$\bigoplus\C_1(H_i)\subseteq\C_1(\bigoplus H_i)$
allow us to define spectra 
$\Max_{\sigma w}\prod\B(H_i)$, $\Max_1\bigoplus\C_1(H_i)$.

The algebra of $n\times n$ complex matrices, which is isomorphic to
$\B(\CC^n)$, is denoted $M_n\CC$. 

\section{Girard couples}

\begin{defin}
A \emph{couple (of quantales)} consists of two quantales $C,Q$
together with a \emph{coupling map} $\phi:C\to Q$ such that $C$ is also
a $Q$-bimodule, $\phi$ is a $Q$-bimodule homomorphism, and
\begin{equation}\label{iota}
\phi(c_1)c_2=c_1\phi(c_2)=c_1c_2
\end{equation}
holds for all $c_1,c_2\in C$.

Assume that $C\stackrel{\phi}{\to}Q$ is a couple.
An element $d\in C$ is said to be \emph{cyclic} if 
$ac\leq d \LRa ca\leq d$ for all $a\in Q,c\in C$. 
The element $d$ is said to be \emph{dualizing} if 
$d\la(a\ra d)=(d\la a)\ra d=a$ for all $a\in Q$ and
$d\la(c\ra d)=(d\la c)\ra d=c$ for all $c\in C$.
In the case where $d$ is both cyclic and dualizing we write
$a\perp c\LRa ac\leq d\LRa ca\leq d$ for $a\in Q,c\in C$ and
$a\comp=a\ra d=d\la a=\bigvee\{ c\in C\mid a\perp c\},
c\comp =c\ra d=d\la c=\bigvee\{ a\in A\mid a\perp c\}$.

A couple $C\stackrel{\phi}{\to}Q$ is said to be:
\emph{strong} if $\phi$ is strong;
\emph{unital} if $Q$ is a unital quantale and $C$ is a unital
$Q$-bimodule; 
\emph{Girard} if it has a cyclic dualizing element.
\end{defin}

\begin{example} \label{examples}
(1) $Q\stackrel{id}{\to}Q$ is clearly a strong couple for any quantale
  $Q$. 
It is unital, or Girard, if and only if $Q$ is unital, or Girard,
  respectively. 

(2) Given an arbitrary unital quantale $Q$, we can construct a Girard
  couple $Q\op\stackrel{0}{\to}Q$ as follows: 
$Q\op$ is equipped with the zero multiplication and the $Q$-bimodule
  structure given by $ac=(a\ra c')',ca=(c'\la a)'$;  
$0$ is the constantly zero map; 
the cyclic dualising element is $e'\in Q\op$.
[Recall that we use $'$ to distinguish elements of $Q$ and $Q\op$.]
Indeed, $a\comp=a\ra e'=\bigvee\{ c\mid ca\leq e'\}=
\bigvee\{ c\mid c'\la a\geq e\}=\bigvee\{ c\mid a\leq c'\}=a'$.
This couple is clearly not strong unless $Q=\{0\}$.

(3) If $C_j\stackrel{\phi_j}{\to}Q_j$ are couples, then so is
$\prod_jC_j\stackrel{(\phi_j)_j}{\to}\prod_jQ_j$. 
Moreover, it is strong, unital, or Girard, if and only if each
component is so. 

(4) Let $R$ be a (unital) ring, $I$ a two-sided ideal, and 
$\Sub R,\Sub I$ the quantales of their additive subgroups. 
Then $\Sub I\subseteq\Sub R$ is a (unital) couple.
\end{example}

\begin{prop} \label{bonus}
Let $C\stackrel{\phi}{\to}Q$ be a couple.
Then

(1) $a(c_1c_2)=(ac_1)c_2$, $(c_1c_2)a=c_1(c_2a)$ and
    $(c_1a)c_2=c_1(ac_2)$ for all $a\in Q,c_1,c_2\in C$. 

(2) $\phi:C\to Q$ is a quantale homomorphism.
\end{prop}
\begin{proof}
(1) Using (\ref{iota}) twice, and the fact that $C$ is a $Q$-bimodule, 
we obtain 
$a(c_1c_2)=a(c_1\phi(c_2))=(ac_1)\phi(c_2)=(ac_1)c_2$.
The proofs of the other two equations are similar.

(2) Using (\ref{iota}) and the fact that $\phi$ is left $Q$-module 
homomorphism, we have
$\phi(c_1c_2)=\phi(\phi(c_1)c_2)=\phi(c_1)\phi(c_2)$. 
\end{proof}

\begin{remark} \label{jeff}
Let $\J$ denote the two-element chain, now regarded not as an object
of $\Sup$ but as a monoidal category in its own right 
(with $\otimes=\wedge$), 
and let $!:0\to 1$ denote the unique non-identity morphism of $\J$. 
Then monoidal functors $F:\J\to\Sup$ are in bijective correspondence
with unital couples of quantales $C\stackrel{\phi}{\to}Q$.  
[By way of comparison, recall that a unital quantale is equivalent to
  a monoid in $\Sup$ which, in turn, is equivalent to a monoidal 
  functor $\T \to \Sup$ where $\T$ is the terminal category.]

The correspondence is given by $C=F_0$, $Q=F_1$, $\phi=F_!$.
The \emph{multiplication natural transformation} of $F$ encompasses all
four binary operations of the couple 
(\eg, its $(1,0)$-component, $F_1 \otimes F_0 \to F_{1\otimes 0}=F_0$,
corresponds to the left action of $Q$ on $C$); 
its naturality is equivalent to the restrictions placed on $\phi$;
the pentagon which its required to satisfy summarises all the 
associativity conditions which a couple satisfies, including those of
Proposition \ref{bonus}(1).   
Similarly, the \emph{unit arrow} of $F$, which must have the form 
$\Two\to F_1$, picks out an element of $Q$;
the triangles which it is required to satisfy assert not only that 
this be a unit for $Q$ but also that it act as a unit on $C$.

A very abstract approach to dualising elements, which can be applied
to a much larger class of monoidal functors, is discussed in a
parallel paper \cite{Eg08}. 
Much of what follows for couples of quantales remains true in the more 
general setting.  
\end{remark}

\begin{prop}
Let $C\stackrel{\phi}{\to}Q$ be a Girard couple.
Then $\phi$ is self-adjoint, \ie\ $\phi^*=\phi$.
\end{prop}
\begin{proof}
The assertion follows from
$c_1\leq\phi\adj(c_2\comp)\LRa\phi(c_1)\leq c_2\comp
\LRa\phi(c_1)c_2=c_1\phi(c_2)\leq d\LRa c_1\leq\phi(c_2)\comp$.
\end{proof}

\begin{remark}
Given a Girard couple $C\stackrel{\phi}{\to}Q$, one can define 
$a\para b=(b\comp a\comp)\comp$ for $a,b\in C \cup Q$.  
The four resultant operations all correspond to the 
\emph{multiplicative join \emph{alias} par} of linear logic.
Collectively, they give $Q\op\stackrel{\phi^*}{\to}C\op$ the structure
of a Girard couple, with neutral element $d'$ and cyclic dualising
element $e'$; 
by the previous proposition, this is isomorphic, as a Girard couple,
to $C\stackrel{\phi}{\to}Q$.  
\end{remark}

\begin{prop}\label{strong}
Let $C\stackrel{\phi}{\to}Q$ be a strong couple of semiunital
quantales. 
Then $\phi$ is an isomorphism on right- and left-sided elements.
\end{prop}
\begin{proof}
We prove the right-sided case.
Let $r\in\R(Q),s\in\R(C)$.
Then $\phi(r1_C)=r\phi(1_C)=r1_Q=r$ and $\phi(s)1_C=s1_C=s$, hence
$\phi|_{\R(C)}$ and $\blank\cdot 1_C$ are mutually inverse
sup-homomorphisms. 
\end{proof}

\begin{prop}\label{strGir}
A Girard couple $C\stackrel{\phi}{\to}Q$ is unital;
if it is also strong, then both $C$ and $Q$ are von Neumann
quantales. 
\end{prop}
\begin{proof}
Let $d\in C$ be a cyclic dualizing element.
All the equalities of \cite[Proposition 6.1.2]{Rose90} can be easily
adapted for Girard couples;  
in particular, $e=d\comp$ is a unit for $Q$.

Now assume that $r\leq d$ for some $r\in R(C)$;
then $r1_C=r1_Q\leq d$, and hence $r\leq 1_Q\comp=0_C$.
That is, the only right- or left-sided element below $d$ is $0$. 
It follows that, for all pairs $r\in\R(C),l\in\LL(Q)$, and for all
pairs $r\in\R(Q),l\in\LL(C)$, $a lr\leq d\LRa lr=0$.
Thus $lr=l\phi(r)=\phi(l)r$ for $r\in\R(C),l\in\LL(C)$ and Proposition
\ref{strong} entail that $C$ is von Neumann.
The previous proposition also entails
$lr=0_Q\LRa\phi\adj(l)\phi\adj(r)=\phi\adj(0_Q)=0_C$
for $r\in\R(Q),l\in\LL(Q)$; 
hence $Q$ is also von Neumann.
\end{proof}

\begin{corol}
Every Girard quantale is von Neumann and the Girard duality extends
the von Neumann duality. 
\end{corol}

\begin{theorem}
Let $S$ be a sup-lattice. 
Then the assignment
$$(x\otimes y')(u\otimes v')=
\begin{cases}
  0 & \text{if } u\leq y, \\
  x \otimes v' & \text{otherwise} 
\end{cases}$$
defines a quantale structure on $S\otimes S\op$ which will be denoted
$\C(S)$. 

The assignment
$$ \phi(x\otimes y')=\rho_x\lambda_y $$
defines a strong Girard couple $\C(S)\stackrel{\phi}{\to}\Q(S)$
with a cyclic dualizing element
$$d=\bigvee_{x\in S}(x\otimes x').$$
\end{theorem}
\begin{proof}
The given binary operation is clearly associative and distributive on
generators of $S\otimes S\op$. 
For example,
\begin{align*}
(x\otimes\bigvee y_i')(u\otimes v') &= (x\otimes(\bigwedge y_i)')(u\otimes v')
&= \begin{cases}
0 & \text{if }\forall i\ u\leq y_i \\
x\otimes v' & \text{otherwise}
\end{cases}\\
&= \bigvee (x\otimes y_i')(u\otimes v').
\end{align*}
Thus, by the definition of $\otimes$, it extends to all elements of
$S\otimes S\op$. 

$\phi$ too is clearly a well-defined sup-homomorphism since it is
\scary{bilinear} on generators. 
It is strong because
$\phi(1_{\C(S)})=\phi(1\otimes 0')=\rho_1\lambda_0=1_{\Q(S)}$. 

$S$ is a left $\Q(S)$-module with action $\alpha x=\alpha(x)$ and
$S\op$ is a right $\Q(S)$-module with action 
$y'\alpha=\alpha\adj(y)'$. 
Therefore $S\otimes S\op$ carries $\Q(S)$-bimodule structure.
The axiom (\ref{iota}) is obtained as follows
\begin{align*}
\phi(x\otimes y')(u\otimes v') &= \rho_x\lambda_y(u)\otimes v' \\
&= \begin{cases}
0 & \text{if } u\leq y,\\
x\otimes v' & \text{otherwise}
\end{cases}\\
&= (x\otimes y')(u\otimes v')
\end{align*}
and symmetrically for $(x\otimes y')\phi(u\otimes v')$.
Since the operations of $\Q(S)$ are given pointwise, 
the remaining axioms of a couple are evident.

The duality $\C(S)\cong\Q(S)\op$ was proven in \cite{JoTi84} and is 
given by
$$(\lambda_x\vee\rho_y)\comp=x\otimes y'.$$
Namely, for $\alpha\in\Q(S),c\in\C(S)$ we have $\alpha\perp c$ when
$x\otimes y'\leq c$ implies that $\alpha\leq\lambda_x\vee\rho_y$.
We will show that $d$ is a cyclic dualizing element of the couple
$\C(S)\stackrel{\phi}{\to}\Q(S)$. 
We can see that
$$(\lambda_x\vee\rho_y)(u\otimes v')\leq d\LRa
u\leq x\text{ and }y\leq v\LRa (u\otimes v')(\lambda_x\vee\rho_y)\leq d$$ 
whenever $(\lambda_x\vee\rho_y)\neq 1,(u\otimes v')\neq 0$ and 
$$1(u\otimes v')\leq d\LRa u\otimes v'=0\LRa (u\otimes v')1\leq d.$$
Hence
$$(x\otimes y')\ra d=d\la (x\otimes y')=\lambda_x\vee\rho_y$$
and
$$(\lambda_x\vee\rho_y)\ra d=x\otimes y'=d\la (\lambda_x\vee\rho_y).$$
Consequently,
$$\alpha\ra d =
\bigvee_{x\in S}(\alpha\adj(x)\otimes x')=
\left(\bigwedge_{x\in S}(\lambda_{\alpha\adj(x)}\vee\rho_x\right)\comp
=\alpha\comp$$
and similarly
$$d\la\alpha =\bigvee_{x\in S}(\alpha(x)\otimes x')=\alpha\comp$$
for every $\alpha\in\Q(S)$.
The inverse duality $\comp:\C(S)\to\Q(S)$ follows directly from a
general property of adjoints $(\bigvee c_i)\ra d=\bigwedge(c_i\ra d)$. 
Thus $d$ is a cyclic dualizing element of
$\C(S)\stackrel{\phi}{\to}\Q(S)$. 
\end{proof}

We remark that the morphism $\phi:\C(S)\to\Q(S)$ is an instance of a
\emph{mix map} \cite{CocSee-mix}. 
G.~N.~Raney \cite{Ra60} proved that this $\phi$ is an isomorphism if
and only if $S$ satisfies \emph{complete distributivity}:  
$$ \bigwedge_{j \in J} \bigvee_{k \in K} \alpha_{jk} 
   = \bigvee_{f \in K^J} \bigwedge_{j \in J} \alpha_{jf(j)}. $$
We obtain the following statement which has already been
mentioned in \cite{KrPa07}.

\begin{corol}
$\Q(S)$ is a Girard quantale if and only if $S$ is completely distributive.
\end{corol}

\begin{theorem}\label{spectra}
Let $H$ be a Hilbert space. 
Then the assignment $\phi(c)=\cl_{\sigma w}(c)$ 
(\ie\ the ultraweak closure) defines a Girard couple
$\Max_1\C_1(H)\stackrel{\phi}{\to}\Max_{\sigma w}\B(H)$
with a cyclic dualizing element 
$$d=\{ C\in\C_1(H)\mid\tr C=0\}.$$
\end{theorem}
\begin{proof}
The basic idea is that $(A,C)\mapsto\tr(AC)=\tr(CA)$ is a bilinear
form on $\B(H)\times\C_1(H)\to\CC$ which is continuous in each
variable \wrt\ to the appropriate topology.

It is known \cite{KaRo97} that the ultraweakly continuous functionals 
on $\B(H)$ are of the form $\tr(C\cdot\blank)$ for some
$C\in\C_1(H)$, and conversely, the $\|\ \|_1$-norm continuous
functionals on $\C_1(H)$ are of the form $\tr(A\cdot\blank)$ for some
$A\in\B(H)$.  
In the spectra, operators correspond to atoms and functionals to
coatoms. 
More precisely, we work with one-dimensional subspaces and kernels of
functionals. 
Every closed subspace (in the topology considered) can then be
obtained as a join of atoms or meet of coatoms and it is known that
the families of atoms and coatoms separate each other.
From this fact it follows that the assignment
$$a\perp c\LRa(\forall A\in a,C\in c\ \tr(AC)=0)$$
admits a duality between $\Max_{\sigma w}\B(H)$ and $\Max_1\C_1(H)$. 
Moreover, the trace is symmetric on operators and thus also on
subspaces, \ie\ $d$ is cyclic and from the duality it follows that $d$
is dualizing. 
We obtain
$$a\comp=\bigwedge\{\ker\tr(A\cdot\blank)\mid A\in a\}.$$

$\Max_1\C_1(H)$ is a $\Max_{\sigma w}\B(H)$-bimodule since $\C_1(H)$
is a two-sided ideal in $\B(H)$ and both multiplications
$\B(H)\times\C_1(H)\to\C_1(H)$, $\C_1(H)\times\B(H)\to\C_1(H)$ are
continuous. 
The ultraweak topology is weaker than the $\|\ \|_1$-norm topology,
thus it defines a closure on $\Max_1\C_1(H)$. 
From continuity it follows again that
$\cl_{\sigma w}(ac)=a\cl_{\sigma w}(c)$ and hence $\phi$ is a bimodule 
homomorphism.
It is strong because $\cl_{\sigma w}\C_1(H)=\B(H)$.
\end{proof}

\begin{corol}
The spectrum $\Max M_n\CC$ is a Girard quantale.
\end{corol}
\begin{proof}
On a finite-dimensional Hilbert space all operators are trace-class
and the norm, $\|\ \|_1$-norm, and ultraweak topologies coincide,
hence $\phi$ is an isomorphism.
\end{proof}

\begin{prop}
Let $H_i$ be a family of Hilbert spaces.
Then
$\Max_1\bigoplus\C_1(H_i)\stackrel{\phi}{\to}\Max_{\sigma w}\prod\B(H_i)$
is a strong Girard couple.

In particular, the spectrum $\Max A$ of a finite-dimensional C*-algebra
$A$ is a Girard quantale.
\end{prop}
\begin{proof}
Since $\prod\B(H_i)\subseteq\B(\bigoplus H_i)$,
$\bigoplus\C_1(H_i)\subseteq\C_1(\bigoplus H_i)$ 
are closed subalgebras,
we can correctly restrict
$\phi:\Max_1\C_1(\bigoplus H_i)\to\Max_{\sigma w}\B(\bigoplus H_i)$
to $\phi|_{\Max_1\bigoplus\C_1(H_i)}:\Max_1\bigoplus\C_1(H_i)\to
\Max_{\sigma w}\prod\B(H_i)$.
Then all calculations are made with respect to the invariant subspaces
$H_i$ and $\{ C\in\bigoplus\C_1(H_i)\mid\tr(C)=0\}$ provides a duality
between elements of $\bigoplus\C_1(H_i)$ and ultraweakly continuous
functionals restricted to $\prod\B(H_i)$. 
It is true that $\cl_{\sigma w}(\bigoplus\C_1(H_i))=\prod\B(H_i)$.
The rest follows from Theorem \ref{spectra}.

Finite-dimensional C*-algebras are of the form 
$\prod_{i=1}^kM_{n_i}\CC=\bigoplus_{i=1}^kM_{n_i}\CC$,
hence the assertion. 
\end{proof}

\begin{theorem}\label{Gir}
Let $C\stackrel{\phi}{\to}Q$ be a Girard couple.
Then $\phi$ factors through a Girard quantale $G$, \ie\ there are
quantale homomorphisms $\gamma:C\to G$ and $\alpha:G\to Q$ such that
$\phi=\alpha\gamma$. 
Moreover, $C\stackrel{\gamma}{\to}G$ is a couple and the $G$-module
actions are given by restricting scalars along $\alpha$: 
\begin{align*}
gc &= \alpha(g)c, & cg &= c\alpha(g)
\end{align*}
for all $g\in G,c\in C$.

If $\phi$ is strong then $\gamma,\alpha$ can be chosen to be strong.
Consequently, $\R(C)\cong\R(G)\cong\R(Q),\LL(C)\cong\LL(G)\cong\LL(Q)$.
\end{theorem}
\begin{proof}
We follow the idea of \cite[Theorem 6.1.3]{Rose90}.
Let $G=\{ (a,c)\in A\times C\mid\phi(c)\leq a\}$ with joins given
componentwise and multiplication 
$(a_1,c_1)(a_2,c_2)=(a_1a_2,a_1c_2\vee c_1a_2)$.
From definition of a couple we easily check
that $A\times C$ is a quantale.
$G$ is clearly closed under joins and from
$\phi(a_1c_2\vee c_1a_2)=a_1\phi(c_2)\vee\phi(c_1)a_2 \leq a_1a_2$ it
follows that it is a strong subquantale of $A\times C$.

Put $\gamma(c)=(\phi(c),c),\alpha(a,c)=a$.
The projection $\alpha$ is evidently a strong quantale homomorphism. 
Actions $c_1(a,c_2)=c_1\gamma(a,C_2)=c_1a$,
$(a,c_1)c_2=\gamma(a,c_1)c_2=ac_2$ define a $G$-bimodule structure on
$C$. 
Then $\gamma(c_1(a,c_2))=\gamma(c_1a)=(\phi(c_1a),c_1a)=
(\phi(c_1)a,c_1a)=\gamma(c_1)(a,c_2)$ because
$\phi(c_1)c_2=c_1c_2\leq ac_2$.
Similarly we can check the other side and hence $\gamma$ is a
$G$-bimodule homomorphism. 
From the properties of $\phi$ it follows that
$C\stackrel{\gamma}{\to}G$ is also a couple. 
The largest element of $G$ is $(1_Q,1_C)$, thus $\gamma$ is strong
whenever $\phi$ is so. 

Finally, $(1,d)$ is a cyclic element of $G$ and $(e,0)$ is a unit. 
Indeed, we have 
$(a_1,c_1)(a_2,c_2)\leq (1,d)\LRa a_1c_2\vee c_1a_2\leq d
\LRa a_1\perp c_2\text{ and } c_1\perp a_2$
which yields the awaited duality $(a,c)\comp=(c\comp,a\comp)$.

The rest follows from Proposition \ref{strong}.
\end{proof}

\begin{corol}
For every sup-lattice $S$ there exists a Girard quantale $\G(S)$ with
$\R\G(S)\cong S,\LL\G(S)\cong S\op$. 
\end{corol}

\begin{remark}
(1) Rosenthal's Girard quantale $Q\times Q\op$ 
(\cite[Theorem 6.1.3]{Rose90}) arises as our $G$ for the zero Girard
  couple $Q\op\stackrel{0}{\to}Q$ of Example \ref{examples}(2). 

(2) The multiplication in $\G(S)$ can be interpreted as a convolution
  product: let $a=(a_0,a_1),b=(b_0,b_1)\in\G(S)$, then
$$ (ab)_i=\bigvee_{j\wedge k\leq i}a_jb_k $$
for both $i\in\{ 0,1\}$.
\end{remark}

\begin{example}
It is possible to meld all the spectra 
$(\Max_p\C_p(H))_{p\in [1,\infty]}$ 
(where $\Max_{\infty}\C_{\infty}H:=\Max_{\sigma w}\B(H)$)
into a single monoidal functor $F:[0,1]\to\Sup$, thus extending the
framework of Remark \ref{jeff}. 

Here $[0,1]$ is regarded as a thin monoidal category with the 
\emph{\L ukasiewicz multiplication} 
$$ i\&_Lj=\max\{ 0, i+j-1\} $$
so that if $p,q,r\in[1,\infty]$ satisfy 
\begin{equation}\label{pq}
\frac{1}{r}=\min\{ 1,\frac{1}{p}+\frac{1}{q}\}
\end{equation}
and if $f$ denotes the bijection $[1,\infty]\to[0,1]$ 
given by $f(p)=1-\frac{1}{p}$, $f^{-1}(i)=\frac{1}{1-i}$,
then $f(p) \&_L f(q) = f(r)$.

The functor $F$ is given by 
$F_i=\Max_{f^{-1}(i)}\C_{f^{-1}(i)}(H)$ and 
$F(i\to j)=\cl_{f^{-1}(j)}$
(\ie\ the $\|\ \|_{f^{-1}(j)}$-norm closure), 
and its multiplication natural transformation by the well-known fact
that $A\in\C_p(H),B\in\C_q(H)$ implies $AB\in\C_r(H)$ where $r$ is
determined by (\ref{pq}). 

Moreover it is possible to construct a single Girard quantale $G$ from
this data, analogous to that constructed in Theorem \ref{Gir}:  
$$ G=\left\{\left.(a_i)\in\prod_{i\in [0,1]}\Max_{f^{-1}(i)}
\C_{f^{-1}(i)}(H)\ \right| \
\cl_{f^{-1}(j)}a_i\subseteq a_j\text{ for }i\leq j\right\} , $$
$$ (ab)_i=\bigvee_{j\&_L k\leq i}\cl_{f^{-1}(i)}(a_jb_k). $$
\end{example}

\begin{oppr}
We have shown that the spectra of the operator algebras $\B(H)$ 
(and their products) together with the spectra of their preduals 
from Girard couples.
It is natural to ask whether our results can be generalized for all  
W*-algebras, including the ideas of the previous example.
If so, it would be useful to describe essential concepts 
of W*-algebras, \eg\ normal morphisms, by means of
the discussed monoidal functors.
\end{oppr}

\bibliography{david}
\bibliographystyle{acm}

\end{document}